\documentclass[final]{siamltex}

\usepackage{amsmath,amsfonts,amssymb,graphics}

\title{Characterization and computation of
$\mathcal{H}_{\infty}$ norms for time-delay systems}

\author{Wim Michiels \and Suat Gumussoy\thanks{Department of Computer Science, Katholieke Universiteit Leuven,
Belgium,  ({\tt
\{Wim.Michiels,Suat.Gumussoy\}@cs.kuleuven.be})} }

\newtheorem{remark}[theorem]{Remark}

\newtheorem{example}[theorem]{Example}
\newtheorem{algorithm}[theorem]{Algorithm}
\newtheorem{assumption}[theorem]{Assumption}

\makeatletter
\def\Ddots{\mathinner{\mkern1mu\raise\p@
\vbox{\kern7\p@\hbox{.}}\mkern2mu
\raise4\p@\hbox{.}\mkern2mu\raise7\p@\hbox{.}\mkern1mu}}
\makeatother

\newcommand{\Hi}{{\cal H}_\infty}
\newcommand{\Lxi}{\mathcal{L}_{\xi}}

\newcommand{\R}{\mathbb{R}}
\newcommand{\C}{\mathbb{C}}

\begin{document}
\maketitle

\begin{abstract}
We consider the characterization and computation of $\Hi$
norms for a class of time-delay systems. It is well known
that in the finite dimensional case the $\Hi$ norm of a
transfer function can be computed using the connections
between the corresponding singular value curves and the
imaginary axis eigenvalues of a Hamiltonian matrix,
leading to the established level set methods. We show a
similar connection between the transfer function of a
time-delay system and the imaginary axis eigenvalues of
an infinite dimensional linear operator $\Lxi$. Based on
this result, we propose a predictor-corrector algorithm
for the computation of the $\Hi$ norm.
In the prediction step, a finite-dimensional
approximation of the problem, induced by a spectral
discretization of the operator $\mathcal{L}_{\xi}$, and
an adaptation of the algorithms for finite-dimensional
systems, allow to obtain an approximation of the $\Hi$
norm of the transfer function of the time-delay system.
In the next step the approximate results are corrected to
the desired accuracy by solving a set of nonlinear
equations which are obtained from the reformulation of
the eigenvalue problem for the linear
infinite-dimensional operator $\Lxi$ as a finite
dimensional nonlinear eigenvalue problem. These equations
can be interpreted as characterizations of peak values in
the singular value plot. The effects of the
discretization in the predictor step are fully
characterized and the choice of the number of
discretization points is discussed. The paper concludes
with a numerical example and the presentation of the
results of extensive benchmarking.

\end{abstract}

\section{Introduction}

In the field of robust control of linear systems
stability and performance criteria are often expressed by
means of $\mathcal{H}_{\infty}$ norms of appropriately
defined transfer functions \cite{zhou}. Therefore, the
availability of robust methods to compute
$\mathcal{H}_{\infty}$ norms is essential in a computer
aided control system design.

In this article we present an approach to compute the
$\mathcal{H}_{\infty}$ norm of the transfer function
\begin{equation}\label{transfer}
G(j\omega)=C\left(j\omega I-A_0-\sum_{i=1}^m A_i
e^{-j\omega\tau_i}\right)^{-1}B+ D e^{-j\omega\tau_0},
\end{equation}
where $A_i\in\mathbb{R}^{n\times n},\ 0\leq i\leq m$,
$B\in\mathbb{R}^{n\times n_u}$,
$C\in\mathbb{R}^{n_y\times n}$,
$D\in\mathbb{R}^{n_y\times n_u}$ are the system matrices,
and the nonnegative numbers
$(\tau_0,\tau_1,\ldots,\tau_m)$ correspond to
time-delays. The $\mathcal{H}_{\infty}$ norm of the
transfer function (\ref{transfer}) is finite if and only
if the zeros of the equation
\[
\det\left(\lambda I-A_0-\sum_{i=1}^m A_i
e^{-\lambda\tau_i}\right)=0
\]
are confined to the open left half complex plane. Under
this condition it can be expressed as
\[
\|G(j\omega)\|_{\mathcal{H}_{\infty}}= \sup_{\omega\geq
0}\sigma_1(G(j\omega)),
\]
where $\sigma_1(\cdot)$ denotes the largest singular
value~\cite{zhou,wimrobust}.

The commonly used methods for computing $\Hi$ norms and
related robustness measures for systems without delay
belong to the class of level set methods. They are based
on the duality between the singular value plot of the
transfer function and the position of the spectrum of an
appropriately defined Hamiltonian matrix with respect to
the imaginary axis, as expressed in the following result
from \cite{boydbala2} (see also \cite{byers}):
\begin{proposition}\label{propintro}
Let  $G_f(j\omega)=C(j\omega I-A)^{-1}B+D$. Let $\xi> 0$
be such that the matrix $D^TD-\xi^2 I$ is non-singular.
For $\omega\geq 0$, the matrix $G_f(j\omega)$ has a
singular value equal to $\xi$ if and only if $j\omega$ is
an eigenvalue of the matrix
\[
L_{\xi}=\left[\begin{array}{cc}  A& -B
(D^TD-\xi^2 I)^{-1}B^T\\
-C^TC+C^TD (D^TD-\xi^2 I)^{-1}D^TC & -A^T
\end{array}\right].
\]
\end{proposition}
From Proposition~\ref{propintro} we get
\[
\|G_f(j\omega)\|_{\Hi}=\inf \left\{\xi> \sigma_1(D):\
L_{\xi} \mathrm{\ has\ no\ imaginary\ axis\
eigenvalues}\right\}.
\]
This result directly leads a bisection algorithm on the
parameter $\xi$ for computing the $\Hi$ norm of
$G_f(j\omega)$, as outlined in \cite{boydbala2}.
Quadratically convergent algorithms based on a search in
a two-parameter space $(\omega,\xi)$ are presented
in~\cite{boydbala,steinbuch}. A similar algorithm for
computing pseudospectral abscissa for systems without
delays is proposed in \cite{overtoncriscross}.

The approach of the paper to compute the $\Hi$ norm of
(\ref{transfer}) builds on a generalization of
Proposition~\ref{propintro} to time-delay systems. Due to
the fact that a time-delay system is inherently
infinite-dimensional \cite{bookwim}, the singular value
curves of (\ref{transfer}) can no longer be related to
the imaginary axis eigenvalues of a matrix but to the
imaginary axis eigenvalues of an infinite-dimensional
linear operator $\mathcal{L}_{\xi}$, as we shall see.
 This leads to a two-step approach for
the computation of the $\Hi$ norm of (\ref{transfer}). In
the first step (the prediction step), an approximation of
the $\Hi$ norm of (\ref{transfer}) is computed based on a
finite-dimensional approximation of the system, induced
by a discretization of the operator $\mathcal{L}_{\xi}$.
Because this operator is a derivative operator on a
function space with nonlocal boundary condition, the
discretization is done using a spectral method
\cite{trefethenspectral}, well established for this type
of operators, see \cite{breda:nonlocal} and the
references therein. In the next step (the correction
step) the approximation of the $\Hi$ norm is improved up
to the desired accuracy with a local method, by solving a
set of nonlinear equations. These are obtained from the
reformulation of the eigenvalue problem for the linear
infinite-dimensional operator $\Lxi$ as a nonlinear
eigenvalue problem of finite dimension.

The proposed method for computing $\Hi$ norms has several
similarities with some existing methods for computing
characteristic roots of time-delay systems, although the
underlying problems are totally different. First, the
characteristic roots solve an infinite-dimensional linear
eigenvalue problem as well as a finite-dimensional
nonlinear eigenvalue problem (induced by the
characteristic equation), see
\cite{bookwim,phdverheyden}. This may also lead to a
two-step approach, where  approximations of the
characteristic roots are obtained by discretizing the
infinite-dimensional linear eigenvalue problem and
solving the resulting matrix eigenvalue problem in the
first place, and the approximate characteristic roots are
corrected subsequently by Newton iterations on the
nonlinear characteristic equation. Such a
predictor-corrector scheme is implemented in the software
package DDE-BIFTOOL \cite{DDE-biftool}. Second, one of
the common approaches to compute characteristic roots
consists of discretizing the infinitesimal generator of
the time-integration operator (solution operator) that
generates the semi-flow of the solutions, see, e.g.\ the
methods proposed in \cite{breda:03,breda}. The
infinitesimal generator is also a derivative operator
with nonlocal boundary conditions, to which a spectral
discretization is employed in \cite{breda}. It forms the
basis for the computation of characteristic roots by the
package TRACE-DDE.

\smallskip

The structure of the article is as follows. In Section
\ref{sectheo} Proposition~\ref{propintro}
 is generalized to transfer functions of
the form (\ref{transfer}). These connections form the
theoretical basis of the paper. In Section~\ref{par1} the
properties of a finite-dimensional approximation based on
a spectral discretization of the infinite-dimensional
operator $\mathcal{L}_{\xi}$ are discussed. In
Section~\ref{paralg} the predictor-corrector algorithm is
described in detail. Section~\ref{secex} is devoted to
the numerical examples. In Section~\ref{seccon} some
concluding remarks are presented.

\subsection*{Notations and assumptions} The notations are
as follows:

\medskip

\begin{tabular}{ll}
  $\C, \R:$ & the field of the complex and real numbers \\
  $\R_+:$ & set of nonnegative real numbers \\
  $A^{*}:$ & complex conjugate transpose of the matrix $A$ \\
  $A^{-T}:$ & transpose of the inverse matrix of $A$ \\
  $\mathcal{D}(.):$ & domain of an operator \\
  $I,I_n$: & identity matrix of appropriate dimensions,
  of dimensions $n\times n$ \\
  $j$: & imaginary identity \\
  $\sigma_i(A):$ & i$^\textrm{th}$ singular value of $A$,\ $\sigma_1(\cdot)\geq\sigma_2(\cdot)\geq \cdots$\\
  $\lambda_i(A):$ & i$^\textrm{th}$ eigenvalue of $A$,\ $|\lambda_1(\cdot)|\geq|\lambda_2(\cdot)|\geq \cdots$\\
  $\Re(u):$ & real part of the complex number $u$ \\
  $\Im(u):$ & imaginary part of the complex number $u$ \\
  $\bar u:$  & complex conjugate of the complex number
  $u$\\
  $|u|$: & modulus of the complex number $u$\\
  $\det(A):$ & determinant of the matrix $A$ \\
  $A\otimes B:$ & Kronecker product of matrices $A$ and
  $B$
\end{tabular}

\medskip

Throughout the paper the following assumption is made:
\begin{assumption}\label{asintro}
\[
\max_{0\leq i\leq m}{\tau_i}=1.
\]
\end{assumption}
Note that Assumption~\ref{asintro} can be taken without
any loss of generality because the variable $\omega$ and
the system matrices in (\ref{transfer}) can always be
re-scaled. It will allow us to significantly simplify the
notations.


\section{Theoretical basis}\label{sectheo}

\subsection{Relations with a Hamiltonian eigenvalue problem}
The following lemma extends Proposition~\ref{propintro}
to time-delay systems:

\begin{lemma}\label{lemeen}
Let $\xi> 0$ be such that the matrix
\[
D_{\xi}:=D^TD-\xi^2 I
\]
is non-singular. For $\omega\geq 0$, the matrix
$G(j\omega)$ has a singular value equal to $\xi$ if and
only if $\lambda=j\omega$ is a solution of the equation
\begin{equation}\label{nonlinear-eigenv}
\det H(\lambda,\xi)=0,
\end{equation}
where
\begin{equation}\label{defh}
H(\lambda,\xi):=\lambda I-M_0-\sum_{i=1}^m \left(M_i
e^{-\lambda\tau_i}+M_{-i}e^{\lambda\tau_i}\right)-\left(N_1e^{-\lambda\tau_0}+
N_{-1}e^{\lambda\tau_0}\right),
\end{equation}
with
\[
\begin{array}{l}
M_{0}=\left[\begin{array}{cc}  A_0& -B
D_{\xi}^{-1}B^T\\
-C^TC+C^TD D_{\xi}^{-1}D^TC & -A_0^T
\end{array}\right],\\
M_i=\left[\begin{array}{cc} A_i &0\\0&0
\end{array}\right],\ \
M_{-i}=\left[\begin{array}{cc} 0 &0\\0&-A_i^T
\end{array}\right],\ \ 1\leq i\leq N,
\\
N_{1}=\left[\begin{array}{cc} 0 & 0\\ 0 & C^TD
D_{\xi}^{-1} B^T
\end{array}\right],\ \
N_{-1}=\left[\begin{array}{cc} -B D_{\xi}^{-1} D^T C& 0\\
0 & 0
\end{array}\right].
\end{array}
\]
\end{lemma}

\noindent\textbf{Proof.\ } The proof is similar to the
proof of Proposition~22 in \cite{genin}. For all
$\omega\in\mathbb{R}$, we have the relation
\begin{multline}\label{vier}
\det H(j\omega,\xi)\det D_{\xi}(j\omega)= \det
(G^*(j\omega)G(j\omega)-\xi^2
I)\\
\det\left(\left[\begin{array}{cc}j\omega I-
A_0-\sum_{i=1}^m A_i e^{-j\omega\tau_i} &0\\
0& j\omega I+A_0^T+\sum_{i=1}^m A_i^T e^{j\omega\tau_i}
\end{array}\right]\right),
\end{multline}
because both left and right hand side can be interpreted
as expressions for the determinant of the 2-by-2 block
matrix
\[
\left[\begin{array}{cc|c} j\omega I- A_0-\sum_{i=1}^m A_i
e^{-j\omega\tau_i} &0& -B\\
C^TC&j\omega I+A_0^T+\sum_{i=1}^m A_i^T e^{j\omega\tau_i}&C^TD\\
\hline D^T C&B^T&D_{\xi}
\end{array}\right]
\]
using Schur complements. Because $D_{\xi}$ is
non-singular and $G$ is stable, we get from (\ref{vier}):
\[
\det (G^*(j\omega)G(j\omega)-\xi^2
I)=0\Leftrightarrow\det H(j\omega,\xi)=0.
\]
This is equivalent to the assertion of the theorem.
 \hfill $\Box$

For a fixed value of $\xi$, the solutions of
(\ref{nonlinear-eigenv}) can be found by solving the
nonlinear eigenvalue problem
\begin{equation}
H(\lambda,\xi)\ v=0,\ \ \lambda\in\mathbb{C},\
v\in\mathbb{C}^{2n},\ v\neq 0. \label{nonlinear-eigenv2}
\end{equation}
This nonlinear eigenvalue problem can be "linearized" to
an infinite-dimensional {linear} eigenvalue problem.
For this, we consider the space
\[
X:=\mathcal{C}\left(\left[-\max_{0\leq i\leq m}\tau_i,\ \max_{0\leq i\leq m}\tau_i\right],\mathbb{C}^{2n}\right)=\mathcal{C}([-1,\ 1],\mathbb{C}^{2n}),
\]
where we have taken into account Assumption~\ref{asintro}.
We let the operator $\mathcal{L}_{\xi}$ on $X$ be
defined by:
\begin{eqnarray}
\mathcal{D}(\mathcal{L}_{\xi}) &=&\left\{\phi\in X:\
\phi^{\prime}\in X,\ \ \phi^{\prime}(0)=M_{0}\phi(0)  +
\sum_{i=1}^m
(M_i\phi(-\tau_i)+M_{-i}\phi(\tau_i) ) \right. \label{defl1}\\
&& \left.\ \hfill \hspace*{6.5cm}
+N_1\phi(-\tau_0) +N_{-1}\phi(\tau_0)
\right\}, \nonumber\\
\mathcal{L}_{\xi}\ \phi&=&\phi^{\prime},\ \ \
\phi\in\mathcal{D}(\mathcal{L}_{\xi}).\hspace*{8.6cm}\label{defl2}
\end{eqnarray}
The eigenvalue problem for this {linear} operator is
defined as
\begin{equation}\label{linear-eigenv2}
(\lambda I-\mathcal{L}_{\xi}) u=0:\
\lambda\in\mathbb{C},\ u\in X,\ u\neq 0.
\end{equation}
The eigenvalues of $\mathcal{L}_{\xi}$ have a one-to-one
correspondence to the eigenvalues of the nonlinear
eigenvalue problem (\ref{nonlinear-eigenv2}):
\begin{proposition}\label{propcor}
Let $H$ be defined by (\ref{defh}). Let $\xi>0$ be such
that $D^TD-\xi^2 I$ is nonsingular. Then we have
\[
\exists v\in\mathbb{C}^{2n},\ v\neq 0:\
H(\lambda,\xi)v=0\Leftrightarrow \exists u\in X,\ u\neq
0:\ (\lambda I-\mathcal{L}_{\xi})u=0.
\]
Furthermore, if $(\lambda,u)$ satisfies
(\ref{linear-eigenv2}), then $u$ has the form
\begin{equation}\label{defu}
u(\theta)=v e^{\lambda\theta},\ \theta\in[-1,\ 1],\
\end{equation}
where $v\in\mathbb{C}^{2n}$ and $(\lambda,v)$ satisfies
(\ref{nonlinear-eigenv2}). Conversely, if $(\lambda,v)$
satisfies (\ref{nonlinear-eigenv2}) then $(\lambda,u)$
satisfied (\ref{linear-eigenv2}) with $u$ given by
(\ref{defu}).
\end{proposition}

\noindent\textbf{Proof.\ } Assume that
$\mathcal{L}_{\xi}\ u=\lambda u$. From (\ref{defl2}) we
get $u(t)=e^{\lambda t} v,\ t\in[-1,\ 1]$, with
$v\in\mathbb{C}^{2n}$. Taking into account the boundary
condition (\ref{defl1}) we get $H(\lambda,\xi)v=0$.
Conversely, if $H(\lambda,\xi)v=0$, then it is readily
verified that $u\equiv ve^{\lambda\theta},\
\theta\in[-1,\ 1]$, belongs to
$\mathcal{D}(\mathcal{L}_{\xi})$ and satisfies
$(\mathcal{L}_{\xi}-\lambda I)u=0$.~$\Box$

By combining Lemma~\ref{lemeen} and
Proposition~\ref{propcor} we arrive at:

\begin{theorem}\label{theoeen}
Let $\xi> 0$ be such that the matrix $D^TD-\xi^2 I $ is
non-singular. For $\omega\geq 0$, the matrix $G(j\omega)$
has a singular value equal to $\xi$ if and only if
$\lambda=j\omega$ is an eigenvalue of the operator
$\mathcal{L}_{\xi}$, defined by (\ref{defl1}) and
(\ref{defl2}).

\end{theorem}

\begin{corollary}\label{colcompute}
\[
\|G(j\omega)\|_{\mathcal{H}_{\infty}}=\inf
\{\xi>\sigma_1(D^TD):\ \mathrm{operator}\
\mathcal{L}_{\xi} \mathrm{\ has\ no\ imaginary\ axis\
eigenvalues}\}.
\]
\end{corollary}

\subsection{Properties of the eigenvalue problem}

\noindent Although the operator $\mathcal{L}_{\xi}$
generally has an infinite number of eigenvalues, the
number of eigenvalues on the imaginary axis is always
finite. This can be concluded from the following results:
\begin{proposition}\label{propenvelope}
All eigenvalues of $\mathcal{L}_{\xi}$ belong to the set
\begin{multline}\nonumber
\Xi:=\left\{\lambda\in\mathbb{C}:\ |\lambda|\leq
\|M_0\|+\sum_{i=1}^m \left(\|M_i\|
e^{-\Re\lambda\tau_i}+\|M_{-i}\|e^{\Re\lambda\tau_i}\right)
\right.
\\
\left. +\|N_1\| e^{-\Re\lambda\tau_0}+\|N_{-1}\|
e^{\Re\lambda\tau_0}\right\}.
\end{multline}
\end{proposition}
\noindent\textbf{Proof.\ } From the identity $\det
H(\lambda,\xi)=0$ we get
\[
|\lambda|\leq \|M_0\|+\sum_{i=1}^m \left(\|M_i\|
\left|e^{-\lambda\tau_i}\right|+\|M_{-i}\|\left|e^{\lambda\tau_i}\right|\right)+\|N_1\|
\left|e^{-\lambda\tau_0}\right|+\|N_{-1}\|
\left|e^{\lambda\tau_0}\right|.
\]
The proposition follows. \hfill $\Box$

\smallskip

Note that $\lambda$ appears in both left and right hand side of the inequality that defines the set $\Xi$. If one considers the elements of $\Xi$ with a given real part, say $\Re(\lambda)=p,\ p\in\mathbb{R}$, then the bound on their modulus $|\lambda|$ is obtained by evaluating the right hand side for $\Re(\lambda)=p$ (a value smaller than $|p|$ implies the set $\Xi$ has no intersection with the line $\Re(\lambda)=p$).

\medskip

\begin{corollary}\label{propfinite}
For all $c>0$, the number of eigenvalues of
$\mathcal{L}_{\xi}$ in the strip
\begin{equation}\label{strip}
\{\lambda\in\mathbb{C}:\ -c<\Re(\lambda)<c\}
\end{equation}
is finite.
\end{corollary}

\noindent\textbf{Proof.\ } Proposition~\ref{propenvelope}
implies that the eigenvalues in the strip (\ref{strip})
can be constrained to  a compact set. Because the
function $H(\cdot,\xi)$ is analytic this number is
finite.\hfill $\Box$

 The set of eigenvalues of $\mathcal{L}_{\xi}$ is
symmetric w.r.t.\ the imaginary axis, as expressed in the
following proposition (for comparison, in the delay-free
case the operator $\mathcal{L}_{\xi}$ reduces to a
Hamiltonian matrix):
\begin{proposition}\label{propsymmetric}
A complex number $\lambda$ is an eigenvalue of
$\mathcal{L}_{\xi}$ if and only if $-\bar\lambda$ is an
eigenvalue of $\mathcal{L}_{\xi}$.
\end{proposition}

\noindent\textbf{Proof.\ } It can be directly verified
that
\[
H(-\bar\lambda,\xi)=-\left( \left(
\left[\begin{array}{rr}0&-1\\1&0
\end{array}\right]\otimes I\right)H(\lambda,\xi)
\left(\left[\begin{array}{rr}0&1\\-1&0
\end{array}\right]\otimes I\right)\right)^*,
\]
hence,
\begin{equation}\label{symmetrie}
\det H(-\bar\lambda,\xi)=\left(\det
H(\lambda,\xi)\right)^*
\end{equation}
and
\[
\det H(-\bar\lambda,\xi)=0 \Leftrightarrow\det
H(\lambda,\xi)=0.
\]
This is equivalent to the statement to be proven. \hfill
$\Box$

\subsection{Interpretation in a two parameter
space}\label{par2par}

We consider the equation
\begin{equation}\label{central}
\det H(j\omega,\xi)=0
\end{equation}
as the central equation, and look at its solutions in the
two parameter space
$(\omega,\xi)\in\mathbb{R}_+\times\mathbb{R}_+$. From
Lemma~\ref{lemeen} and Proposition~\ref{propcor} we have
the following corollary:

\begin{theorem} Let $\omega\geq 0$ and $\xi\geq 0$ such
that the matrix $D_{\xi}$ is nonsingular. The following
statements are equivalent:
\begin{enumerate}
\item $\det H(j\omega,\xi)=0$;
\item $j\omega$ is an eigenvalue of the operator $\mathcal{L}_{\xi}$;
\item $\xi$ is a singular value of the matrix $G(j\omega)$.
\end{enumerate}\label{corgraph}
\end{theorem}


\smallskip

\begin{figure}
\begin{center}
\resizebox{11cm}{!}{\includegraphics{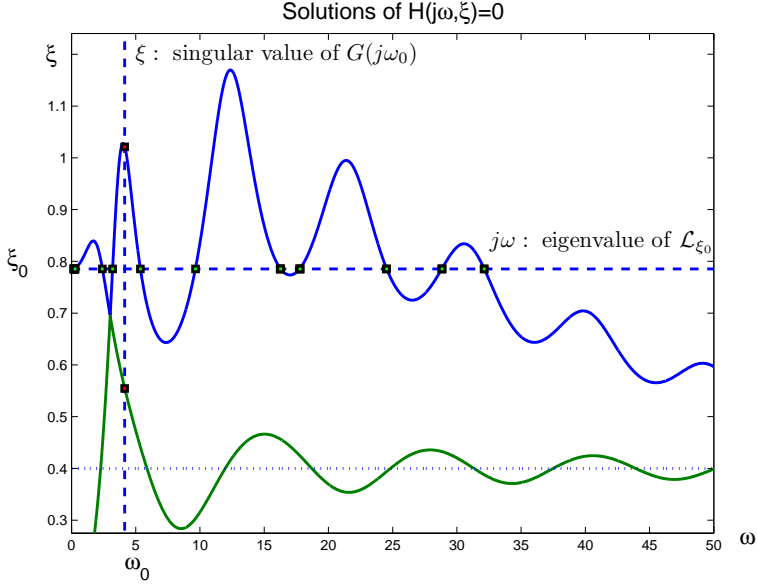}}
 \caption{\label{figA} Solutions of the equation
 (\ref{central}), for the problem data (\ref{example1}).}
\end{center}
\end{figure}

Graphically, the equation (\ref{central}) defines a set
of curves in the $(\omega,\xi)$ parameter space. The
intersections with horizontal lines ($\xi$ fixed) can be
found by computing the imaginary axis eigenvalues of
$\mathcal{L}_{\xi}$. Similarly, the intersections with
vertical lines ($\omega$ fixed) can be found by computing
the singular values of $G(j\omega)$. This is shown in
Figure~\ref{figA}.

\section{Finite-dimensional approximation}\label{par1}

The numerical methods for computing
$\mathcal{H}_{\infty}$ norms presented in
Section~\ref{paralg} are strongly based on
Corollary~\ref{colcompute}. Because the operator
$\mathcal{L}_{\xi}$, defined by
(\ref{defl1})-(\ref{defl2}), is infinite-dimensional,
these algorithms will involve a discretization of this
operator. In this section we outline a discretization
approach and discuss its properties.

\subsection{Discretization}\label{seccdisc}
Following the approach of \cite{breda,breda:nonlocal}, we
discretize the operator $\mathcal{L}_{\xi}$  using a
\emph{spectral method} (see, e.g.
\cite{trefethenspectral}).

Given a positive integer $N$, we consider a mesh
$\Omega_N$ of $2 N+1$ distinct points in the interval
$[-1,\ 1]$:
\begin{equation}\label{defmesh}
\Omega_N=\left\{\theta_{N,i},\ i=-N,\ldots,N\right\},
\end{equation}
where
\[
-1\leq\theta_{N,-N}<\ldots<
\theta_{N,-1}<\theta_{N,0}=0<\theta_{N,1}<\cdots<\theta_{N,N}\leq1
\]
and
\begin{equation}\label{assymmetric}
\theta_{N,-i}=-\theta_{N,i},\ \  i=1,\ldots,N.
\end{equation}

The mesh $\Omega_N$ allows us to replace the continuous
space $X$ with a space $X_N$ of discrete functions. More
precisely, a function $\phi\in X$ is discretized into a
block vector $x=[x_{-N}^T\cdots\ x_{N}^T ]^T\in X_N$ with
components
\[
x_i=\phi(\theta_{N,i})\in\mathbb{C}^{2n},\ \
i=-N,\ldots,N.
\]
When defining $\mathcal{P}_N x,\ x\in X_N$ as the unique
$\mathbb{C}^{2n}$ valued interpolating polynomial of
degree less than or equal to $ 2N$ satisfying
\[
\mathcal{P}_N x (\theta_{N,i})=x_{i},\ \ i=-N,\ldots,N,
\]
we can approximate the operator $\mathcal{L}_{\xi}$ over
$X$ with the matrix $\mathcal{L}_{\xi}^N:\ X_N\rightarrow
X_N$, defined as
\begin{equation}\label{defldisc}
\begin{array}{lll}
\left(\mathcal{L}_{\xi}^N\
x\right)_i=&\left(\mathcal{P}_N
x\right)^{\prime}(\theta_{N,i}), & i=-N,\ldots,-1,\\
\left(\mathcal{L}_{\xi}^N\ x\right)_0=
&M_0 \mathcal{P}_N x(0)+\sum_{i=1}^m \left(M_i
\mathcal{P}_N x(-\tau_i)+ M_{-i} \mathcal{P}_N
x(\tau_i)\right)
\\
& +N_1 \mathcal{P}_N x(-\tau_0)+N_{-1} \mathcal{P}_N
x(\tau_0), &
\\
\left(\mathcal{L}_{\xi}^N\
x\right)_i=&\left(\mathcal{P}_N
x\right)^{\prime}(\theta_{N,i}), & i=1,\ldots,N.\\
\end{array}
\end{equation}
An \emph{explicit expression} for the elements of the matrix $\mathcal{L_{\xi}^N}$ can be obtained by using the Lagrange representation of $\mathcal{P}_N x$,
\begin{equation}\label{lagrange}
\begin{array}{l}
\mathcal{P}_N x=\sum_{k=-N}^N l_{N,k}\ x_k,\ \ \
\end{array}
\end{equation}
where the Lagrange polynomials $l_{N,k}$ are real valued
polynomials of degree $2N$ satisfying
\[
l_{N,k}(\theta_{N,i})=\left\{\begin{array}{ll}1 & i=k,\\
0 & i\neq k.
\end{array}\right.
\]
By substituting (\ref{lagrange}) in (\ref{defldisc}) we obtain the expression
\[
\mathcal{L}_{\xi}^N=
\left[\begin{array}{lll}
d_{-N,-N} &\hdots & d_{-N,N} \\
\vdots & & \vdots \\
d_{-1,-N} &\hdots & d_{-1,N} \\
a_{-N} & \hdots & a_N\\
d_{1,-N} &\hdots & d_{1,N} \\
\vdots & & \vdots \\
d_{N,-N} &\hdots & d_{N,N}
\end{array}\right]\in\mathbb{R}^{(2N+1)(2n)\times(2N+1)2n},
\]
where
\[
\begin{array}{lll}
d_{i,k}&=&l^{\prime}_{N,k}(\theta_{N,i}) I,\ \ \ \
i\in\{-N,\ldots,-1,1,\ldots,N\},\ k\in\{-N,\ldots,N\}, \\
a_0&=& M_0+\sum_{k=1}^m \left(M_k
l_{N,0}(-\tau_k)+M_{-k}l_{N,0}(\tau_k)\right)+N_1
l_{N,0}(-\tau_0)+N_{-1}l_{N,0}(\tau_0), \\
a_{i}&=&\sum_{k=1}^m \left(M_k
l_{N,i}(-\tau_k)+M_{-k}l_{N,i}(\tau_k)\right)+N_1
l_{N,i}(-\tau_0)+N_{-1}l_{N,i}(\tau_0),\\
&&i\in\{-N,\ldots,-1,1,\ldots,N\}.
\end{array}
\]
It is important to note that all the problem specific
information and the parameter $\xi$ are concentrated in
the middle row of $\mathcal{L}_{\xi}^N$, i.e.\ the
elements $(a_{-N},\ldots,a_N)$, while all other elements
of $\mathcal{L}_{\xi}^N$ can be computed beforehand.

We outline some properties of the matrix
$\mathcal{L}_{\xi}^N$.
First, analogously to the continuous case the (linear)
eigenvalue problem for $\mathcal{L}^N_{\xi}$,
\begin{equation}\label{infd}
\mathcal{L}_{\xi}^N\ x=\lambda x,\ \lambda\in\mathbb{C},\
x\in\mathbb{C}^{(2N+1)2n},\ x\neq 0,
\end{equation}
has a nonlinear eigenvalue problem of dimension $2n$ as
counterpart. To clarify this we need the following
definition:
\begin{definition}\label{Defpn}
For $\lambda\in\mathbb{C}$, let $p_N(\cdot;\ \lambda)$ be
the polynomial of degree $2N$ satisfying
\begin{equation}\label{defpn}
\begin{array}{l}
p_N(0;\ \lambda)=1,\\
p_N^{\prime}(\theta_{N,i};\ \lambda)=\lambda
p_N(\theta_{N,i};\ \lambda),\ \
i\in\{-N,\ldots,-1\}\cup\{1,\ldots,N\}.
\end{array}
\end{equation}
\end{definition}
Note that the polynomial $p_N(t;\ \lambda)$ is an
approximation of $\exp(\lambda t)$ on the interval $[-1;\
1]$. Indeed, the first equation of (\ref{defpn}) is an
interpolation requirement at zero, the other equations
are collocation conditions for the differential equation
$\dot z=\lambda z$, of which $\exp(\lambda t)$ is a
solution. We can now state:
\begin{proposition}\label{propnl}
The following statements are equivalent:
\[
\exists x\in \mathbb{C}^{(2N+1)2n},\ x\neq 0:\
\left(\lambda I-\mathcal{L}^N_{\xi}\right)x=0
\Leftrightarrow \exists v\in\mathbb{C}^{2n},\ v\neq 0:\
H_N(\lambda,\xi)v=0 ,
\]
where
\begin{equation}\label{defhn}
\begin{array}{llr}
H_N(\lambda,\xi)&:=&\lambda I-M_0-\sum_{i=1}^m \left(M_i
p_N(-\tau_i;\ \lambda)+M_{-i}p_N(\tau_i,\
\lambda)\right)\ \
\\
&&-\left(N_1 p_N(-\tau_0;\ \lambda)+N_{-1}p_N(\tau_0;\
\lambda)\right).
\end{array}
\end{equation}
%
\end{proposition}

\noindent\textbf{Proof.\ }
Using (\ref{defldisc}) the expression (\ref{infd}) can be
written as:
\begin{eqnarray}
(\mathcal{P}_Nx)^{\prime}(\theta_{N,i})=\lambda
x_i=\lambda \mathcal{P}_Nx(\theta_{N,i}),\  i\in\{-N,\ldots,-1,1,\ldots,N\}, \label{een} \\
M_0 \mathcal{P}_N x(0)+\sum_{i=1}^m \left(M_i
\mathcal{P}_N x(-\tau_i)+ M_{-i} \mathcal{P}_N
x(\tau_i)\right)+N_1 \mathcal{P}_N x(-\tau_0)\nonumber\\
+N_{-1} \mathcal{P}_N x(\tau_0)=\lambda
x_0=\lambda\mathcal{P}_N x(0) . \label{twee}
\end{eqnarray}
From  $\mathcal{P}_N x(0)=x_0$ and (\ref{een}) it follows
that
\begin{equation}\label{collocation}
\mathcal{P}_N x(\cdot)=p_N(\cdot;\ \lambda) x_0.
\end{equation}
When substituting (\ref{collocation}) in (\ref{twee}) we
arrive at $H_N(\lambda,\xi)\ x_0=0$. \hfill$\Box$

Notice that $H_N(\lambda,\xi)$ can be obtained from
$H(\lambda,\xi)$, by making the substitution
\[
e^{-\lambda\tau_i}\leftarrow p_N(-\tau_i;\ \lambda),\
i=-N,\ldots,N.
\]
As we shall see in \S\ref{parpredict} the functions
$p_N(-\tau_i;\ \lambda)$ are proper rational functions of
the parameter $\lambda$. Thus, the effect of a spectral
discretization of the operator $\mathcal{L}_{\xi}$ can be
interpreted as the effect of a \emph{rational
approximation} of the exponential functions in
$H(\lambda,\xi)$.

Second, the spectral property
described in Proposition~\ref{propsymmetric} is
preserved, due to the symmetry of the grid:
\begin{proposition}\label{propsym}
A complex number $\lambda$ is an eigenvalue of
$\mathcal{L}^N_{\xi}$ if and only if $-\bar\lambda$ is an
eigenvalue of $\mathcal{L}^N_{\xi}$.
\end{proposition}

\noindent\textbf{Proof.\ }
The property (\ref{assymmetric}) of the grid assures that
\[
p_N(-\tau_i;\ \lambda)=p_N(\tau_i;\ -\lambda),\ \forall\lambda\in\mathbb{C},\ \forall i\in\{0,\ldots,m\}.
\]
Using this result, the similar arguments as in the proof of
Proposition~\ref{propsymmetric} lead us to
\[
\det H_N(-\bar\lambda,\xi)=0\Leftrightarrow \det H_N(\lambda,\xi)=0,
\]
which is equivalent to the statement of the proposition.
 \hfill $\Box$

\subsection{Interpretation in a two parameter space}
\label{twopardis}

Similarly as in \S\ref{par2par} we characterize the
solutions of the two-parameter problem
\begin{equation}\label{twopardisc}
\det H_N(j\omega,\xi)=0,\ \ \
\end{equation}
where $\omega\geq 0$ and $\xi\geq 0$.
 The counterpart of Theorem~\ref{corgraph} reads
as:
\begin{theorem}\label{proptwopard}
Let $\omega\geq 0$ and let $\xi\geq 0$ such that the
matrix $D_{\xi}$ is nonsingular. The following statements
are equivalent:
\begin{enumerate}
\item $\det H_N(j\omega,\xi)=0$;
\item $j\omega$ is an eigenvalue of $\mathcal{L}_{\xi}^N$;
\item $\xi^2$ is an eigenvalue of the matrix
\begin{equation}\label{defMN}
{\small
\begin{array}{l}
M_N(j\omega):=\\
\left[\begin{array}{cc}
X_N(j\omega)
+X_N(j\omega)^*+Y_N(j\omega)^*Y_N(j\omega)+D^TD &
\sqrt{r(\omega)}\ X_N(j\omega)^* \\
\sqrt{r(\omega)}\ X_N(j\omega) & D^TD
\end{array}\right]
\end{array}
}
\end{equation}
where
\[
\begin{array}{l}
X_N=D^T C \left(j\omega I-A_0-\sum_{i=1}^m A_i
p_N(-\tau_i;\ j\omega)\right)^{-1}B p_N(\tau_0;\
j\omega),
\\
Y_N=C\left(j\omega I-A_0-\sum_{i=1}^m A_i
p_N(-\tau_i;\ j\omega)\right)^{-1}B,\\
r_N(\omega)=\frac{1}{|p_N(-\tau_0;\ j\omega)|^2}-1.
\end{array}
\]
\end{enumerate}
\end{theorem}

\noindent\textbf{Proof.\ }
The equivalence between the first and second statement
corresponds to Proposition~\ref{propnl}. Therefore, it is
sufficient to prove the equivalence between the first and
the third statement.

We can express $H_N(j\omega,\xi)$ as
{\small
\[\begin{array}{r}
 H_N(j\omega,\xi)=
\left[\begin{array}{cc}
j\omega I-A_0-\sum_{i=1}^m A_i\ p_N(-\tau_i; j\omega) &
0\\
C^TC & j\omega I+A^T+\sum_{i=1}^m A_i^T\ p_N(\tau_i;
j\omega)
\end{array}\right]\ \ \\
+\left[\begin{array}{cc}
B & 0\\
0 & -C^T D p_N(-\tau_0;\ j\omega)
\end{array}\right]
\left[\begin{array}{cc}
(D^T D-\xi^2 I)^{-1} & 0\\
0 & (D^T D-\xi^2 I)^{-1}
\end{array}\right]\ \
\\
\cdot\left[\begin{array}{cc}
p_N(\tau_0;\ j\omega) D^TC & B^T\\
\frac{1}{p_N(-\tau_0;\ j\omega)} D^TC & B^T
\end{array}\right].
\end{array}\]
}
Since the transfer function (\ref{transfer}) has no
poles on the imaginary axis, we have
\[
\det H_N(j\omega,\xi)=0 \Leftrightarrow \det \tilde
H_N(j\omega,\xi)=0,
\]
where {\small
\begin{multline}\nonumber
\tilde H_N(j\omega,\xi):= \xi^2 I-
\left[\begin{array}{cc} D^TD &0 \\ 0 &
D^TD\end{array}\right]
-\left[\begin{array}{cc}
p_N(\tau_0;\ j\omega) D^TC & B^T\\
\frac{1}{p_N(-\tau_0;\ j\omega)} D^TC & B^T
\end{array}\right] \cdot \\
\left[\begin{array}{cc}
j\omega I-A_0-\sum_{i=1}^m A_i\ p_N(-\tau_i; j\omega) &
0\\
C^TC & j\omega I+A^T+\sum_{i=1}^m A_i^T\ p_N(\tau_i;
j\omega)
\end{array}\right]^{-1}\\
\cdot \left[\begin{array}{cc}
B & 0\\
0 & -C^T D p_N(-\tau_0;\ j\omega)
\end{array}\right].
\end{multline}

} Using an explicit formula for the inverse of a
two-by-two block matrix, we obtain:
{\small
\[
\begin{array}{lll}
\tilde H_N(j\omega,\xi)&=& \xi^2 I-
\left[\begin{array}{ll} D^TD+X_N(j\omega)+Y_N(j\omega)^*
Y_N(j\omega) & X_N(j\omega)^* \\
\frac{1}{|p_N(-\tau_0;\ j\omega)|^2}
X_N(j\omega)+Y_N(j\omega)^* Y_N(j\omega) &D^TD+
X_N(j\omega)^*
\end{array}\right].
\end{array}
\]
} By elementary row and column operations we arrive at
\[
\det\tilde H_N(j\omega,\xi)=\det\left(\xi^2
I-M_N(j\omega)\right),
\]
hence,
\[
\det H_N(j\omega,\xi)=0\Leftrightarrow \det\left(\xi^2
I-M_N(j\omega)\right)=0.
\]
This completes the proof.\hfill $\Box$

\begin{figure}
\begin{center}
\resizebox{11cm}{!}{\includegraphics{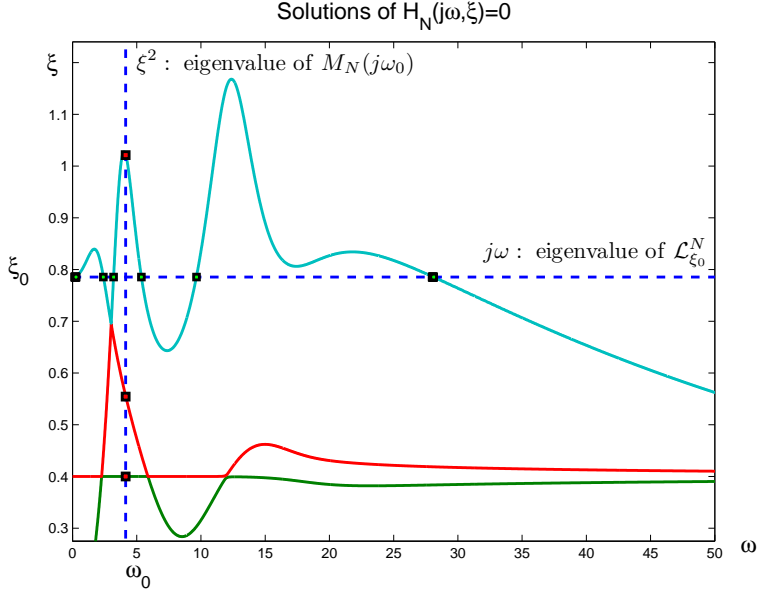}}
 \caption{\label{figB} Solutions of the equation (\ref{twopardisc}) with $N$=10,
 for the problem data (\ref{example1}).}
\end{center}
\end{figure}

The graphical interpretation of Theorem~\ref{proptwopard}
is as follows. The intersections between the curves in
the $(\omega,\xi)$ plane, defined by (\ref{twopardisc}),
with horizontal lines ($\xi$ fixed) can be found by
computing the imaginary axis eigenvalues of the matrix
$\mathcal{L}^N_{\xi}$. Similarly, the intersections with
vertical lines ($\omega$ fixed) can be found by computing
the positive real eigenvalues of the matrix
$M_N(j\omega)$. This is illustrated in Figure~\ref{figB}.

\medskip

In the remainder of the section, we discuss some spectral
properties of the matrix $M_N(j\omega)$, defined in
(\ref{defMN}).

\begin{proposition}\label{propMN}
The following statements hold.
\begin{enumerate}
\item The eigenvalues of $M_N(j\omega)$ are real or appear in
complex conjugate pairs.
\item The matrix $M_N(j\omega)$ has $2n_u$ eigenvalues
$\lambda_k(M_N(j\omega),\ 1\leq k\leq 2n_u$, satisfying
\[
\lim_{\omega\rightarrow\infty}
\lambda_k(M_N(j\omega))=(\sigma_l(D))^2\, \mathrm{\ for\
some\ } l\in\{1,\ldots,n_u\}.
\]
\item If $r(\omega)\geq 0$, or, equivalently, $|p_N(-\tau_0,\ j\omega)|\leq 1$,  then
the matrix $M_N(j\omega)$ is Hermitian.
\item If $\tau_0=0$ then
\[
M_N(j\omega)= \left[\begin{array}{cc}
G_N(j\omega)^*G_N(j\omega) & 0\\
0 & D^T D
\end{array}\right],
\]
where
\begin{equation}\label{defGN}
G_N(j\omega):=C\left(j\omega I-A_0-\sum_{i=1}^m A_i
p_N(-\tau_i;\ j\omega)\right)^{-1}B+D,
\end{equation}
with $p_N$  defined by (\ref{defpn}).
\end{enumerate}
\end{proposition}

\noindent\textbf{Proof.\ }
 The first assertion follows from
\[
M_N(j\omega)= \left[\begin{array}{rc}
\mathrm{sign}(r(\omega)) I& 0 \\
0 & I
\end{array}\right]
M_N(j\omega)^* \left[\begin{array}{rc}
\mathrm{sign}(r(\omega)) I& 0 \\
0 & I
\end{array}\right].
\]
The assertion is implied by
\[
\lim_{\omega\rightarrow\infty} M_N(j\omega)=
\left[\begin{array}{cc} D^TD &0 \\
0 & D^TD
\end{array}\right].
\]
The third assertion is trivial. The fourth assertion is
due to $p_N(-\tau_0;\ j\omega)\equiv 1$ and
$r(\omega)\equiv 0$ if $\tau_0=0$. \hfill $\Box$

The fourth assertion of Proposition~\ref{propMN} is of
particular interest, because it shows that under the
condition $\tau_0= 0$ the curves defined by
(\ref{twopardisc}) can be interpreted as the singular
value curves of a rational approximation $G_N(j\omega)$
of $G(j\omega)$.
%
%
%
%
%

To conclude the section we summarize the relations
between the eigenvalue problems defined in
Sections~\ref{sectheo}-\ref{par1} in
Table~\ref{tabelleke}.

\begin{table}[h]
{\scriptsize{
\begin{tabular}{|l|l||l||l|l|}\hline
&implicit two & horizontal search &vertical search & corresponding \\
&parameter problem &(explicit in $\omega$) &(explicit in
$\xi$)
&$\mathcal{H}_{\infty}$-problem\\
\hline\hline
continuous& $\det H(j\omega,\xi)=0$
 &
 $(j\omega I-\mathcal{L}_{\xi}) u=0$
&   $\left(\xi^2 I-G(j\omega)^* G(j\omega)\right)v=0$ &
$\|G(j\omega)\|_{\mathcal{H}_{\infty}}$\\
(\S\ref{par2par})&&&&
\\
\hline discretized & $\det H_N(j\omega,\xi)=0$
 &$(j\omega I-\mathcal{L}^N_{\xi})x=0$ &
 $(\xi^2 I-M_N(j\omega))v=0$

 &
$\|G_N(j\omega)\|_{\mathcal{H}_{\infty}}$ \\
(\S\ref{twopardis}) & &&&if $\tau_0=0$
\\ \hline
\end{tabular}
} }
\caption{\label{tabelleke} Relations between
$\mathcal{L}_{\xi}$, $H(j\omega,\xi)$ and
$\|G(j\omega)\|_{\mathcal{H}_{\infty}}$, as well as their
discrete counterparts. The latter are all induced by a
spectral discretization of the operator
$\mathcal{L}_{\xi}$ into the matrix
$\mathcal{L}^N_{\xi}$.}
\end{table}

\section{Algorithm}\label{paralg}

Similar to the algorithm for the computation of
characteristic roots of time-delay systems implemented in
the package DDE-BIFTOOL \cite{DDE-biftool}, we propose an
algorithm for computing $\mathcal{H}_{\infty}$ norms that
relies on a two-step approach. In the first step,
outlined in \S\ref{parpredict}, we approximate (predict)
the $\mathcal{H}_{\infty}$ norm based on the
discretization of the operator $\mathcal{L}_{\xi}$ into
the matrix $\mathcal{L}_{\xi}^N$. In the second step,
outlined in \S\ref{parcorrect}, we correct the results
based on the reformulation of the eigenvalue problem for
$\mathcal{L}_{\xi}$ as a nonlinear eigenvalue problem of
finite dimension. In \S\ref{parsteplength} we discuss the
choice of the discretization stepsize in the prediction
step.

Throughout this section we assume that the grid
$\Omega_N$, employed in the discretization of
$\mathcal{L}_{\xi}$, consists of Chebyshev extremal
points in the discretization, that is,
\begin{equation}\label{chebpoints}
\theta_{N,i}=\cos\left(\frac{(N-i)\pi}{2N}\right),\ \
i=-N,\ldots,N.
\end{equation}

\subsection{Prediction of the $\mathcal{H}_{\infty}$ norm}
\label{parpredict}

 Inspired by Corollary~\ref{colcompute} a natural
way to approximate
$\|G(j\omega)\|_{\mathcal{H}_{\infty}}$ consists of
computing
\begin{equation}\label{gmax}
g_{\max}(N):=\inf\left\{\xi>\sigma_1(D^TD):\
\mathrm{matrix\ }\mathcal{L}_{\xi}^N \mathrm{\ has\ no\
imaginary\ axis\ eigenvalues}\right\},
\end{equation}
where $N$ is fixed. The following properties lay the
basis for the corresponding algorithms.
\begin{proposition} If the condition
\begin{equation}\label{proppn}
|p_N(-\tau_0;\ j\omega)|\leq 1,\ \forall\omega\geq 0.
\end{equation}
is satisfied, then the following statements hold.
\begin{enumerate}
\item The quantity $g_{\max}(N)$, defined in (\ref{gmax}), is
finite.
\item  The matrix $\mathcal{L}_{\xi}^N$ has eigenvalues on the imaginary
axis for all
\[
\xi\in\left(\sigma_1(D),\ g_{\max}(N)\right]
\]
and no eigenvalues on the imaginary axis for
$\xi>g_{\max}(N)$.
\item The matrix $M_N(j\omega)$ has $2n_u$ real eigenvalues
$\lambda_k(M_N(j\omega),\ 1\leq k\leq 2n_u$, for all
$\omega\geq 0$, satisfying
\[
\lim_{\omega\rightarrow\infty}
\lambda_k(M_N(j\omega))=(\sigma_l(D))^2\, \mathrm{\ for\
some\ } l\in\{1,\ldots,n_u\}.
\]
\end{enumerate}
Moreover, if $\tau_0=0$, then
\[
g_{\max}(N)=\|G_N(j\omega)\|_{\mathcal{H}_{\infty}},
\]
where $G_N(j\omega)$ is defined in (\ref{defGN}).
\label{propbasbis}
\end{proposition}

\noindent\textbf{Proof.\ } The proof directly follows
from Proposition~\ref{propMN} and
Theorem~\ref{proptwopard}. \hfill $\Box$

Notice that the condition (\ref{proppn}) is always
satisfied if $\tau_0=0$, which implies $p_N(-\tau_0;\
j\omega)\equiv 1$. If $\tau_0\neq 0$, then it is also
satisfied, by the choice (\ref{chebpoints}) of the grid
$\Omega_N$.

\medskip

The definition (\ref{gmax}) and the properties 1.-2.
described in Proposition~\ref{propbasbis} naturally lead
to a bisection algorithm on the parameter $\xi$ to
compute $g_{\max}(N)$, similar to the algorithm presented
in~\cite{byers}. However, based on the interpretations described in
\S\ref{twopardis} and property 3. of
Proposition~\ref{propbasbis}, the efficiency can be
improved by performing a criss-cross search in the two
parameter space $(\omega,\xi)$, instead of a search in
the parameter $\xi$ only. More precisely, an adaptation
of the algorithm presented in \cite{steinbuch} results in
(see also \cite{overtoncriscross}):

\begin{algorithm}\label{algmiddle} \ \

{\small
 {\em
 \noindent Input: system data, $N$, symmetric
grid $\Omega_N$ defined by (\ref{chebpoints}), candidate
critical frequency $\omega_t$ if available,\\
\hspace*{1cm}tolerance tol\\
 Output: $g_{\max}(N)$
 }
\begin{enumerate}
\item compute a lower bound $\xi_l$ on
$g_{\max}(N)$,\\
e.g.\
$\xi_l=\max\left\{\sigma_1(G(0)),\sigma_1(D),\mathrm{tol}
,
\sqrt{\max(\lambda_1(M_N(j\omega_t)),0)}\right\}$
\item repeat until break
  \begin{itemize}
  \item[2.1] set  $\xi:=\xi_l(1+2\ \mathrm{tol})$
  \item[2.2] compute the set
  of eigenvalues $\mathcal{E}_{\xi}$ of the matrix $\mathcal{L}_{\xi}^N$ on the positive imaginary
  axis,\\
   $\mathcal{E}_{\xi}:=\left\{j\omega^{(1)},j\omega^{(2)},\ldots
   \right\}$, with $0\leq \omega^{(1)}<\omega^{(2)}<\ldots$
  \item[2.3] if $\mathcal{E}_{\xi}=\phi$, break\\
  else\\
  $\ \ \ \mu^{(i)}:=\sqrt{\omega^{(i)}\omega^{(i+1)}},\
  i=1,2,\ldots$\\
%
%
  \\
  set
  $\xi_l:= \max_i\sqrt{\max\left(\lambda_1(M_N(j\mu^{(i)})),\sigma_1(D)^2\right)}$,
  \end{itemize}
 \item[] \hspace*{-0.4cm}\{result: estimate $(\xi+\xi_l)/2$
 for $g_{\max}(N)$\}
\end{enumerate}
}
\end{algorithm}

\begin{figure}
\begin{center}
\resizebox{10cm}{!}{\includegraphics{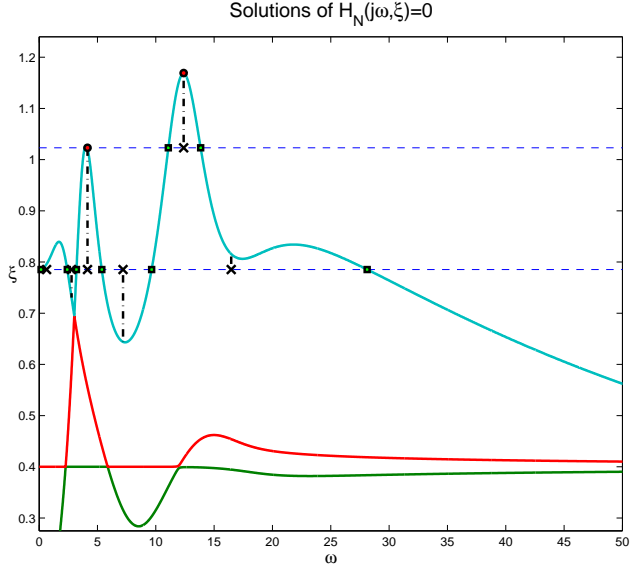}}
 \caption{\label{figC} Some steps of Algorithm~\ref{algmiddle},
 for the problem data (\ref{example1}).}
\end{center}
\end{figure}

The underlying idea is illustrated in Figure~\ref{figC},
where some steps of the algorithm are visualized.


\smallskip

Algorithm~\ref{algmiddle} relies on the evaluation of the
matrix $M_N(\lambda)$, and, hence, on the evaluation of
the functions
\[
\lambda\mapsto p_N(\pm\tau_i;\ \lambda),\ i=0,\ldots,m
\]
for {specific values} of $\lambda$. For this, we
represent $p_N(\cdot;\ \lambda)$ in a polynomial basis:
\[
p_N(t;\ \lambda)=\sum_{i=0}^{2N}\alpha_i
T_i\left(t\right),\ \ \ t\in[-1,\ 1],
\]
where we suppress the dependence of the coefficients
$\alpha_i$ on $\lambda$ in the notation. For $\lambda\neq
0$, the conditions (\ref{defpn}) can be written as
\begin{equation}\label{solvepn}
\left(\lambda T-U\right) \left[\begin{array}{c}
\alpha_0\\
\vdots\\
\alpha_{2N}
\end{array}\right] =R,
\end{equation}
where
{\scriptsize
\[
T=\left[\begin{array}{lll}
T_0(\theta_{N,-N}) &\cdots & T_{2N}(\theta_{N,-N})\\
\vdots & & \vdots\\
T_0(\theta_{N,-1}) &\cdots & T_{2N}(\theta_{N,-1})\\
T_0(0) &\cdots & T_{2N}(0)\\
T_0(\theta_{N,1}) &\cdots & T_{2N}(\theta_{N,1})\\
\vdots & & \vdots\\
T_0(\theta_{N,N}) &\cdots & T_{2N}(\theta_{N,N})\\
\end{array}\right],\ \ U=
\left[\begin{array}{lll}
T_0'(\theta_{N,-N}) &\cdots & T_{2N}'(\theta_{N,-N})\\
\vdots & & \vdots\\
T_0'(\theta_{N,-1}) &\cdots & T_{2N}'(\theta_{N,-1})\\
0 &\cdots & 0\\
T_0'(\theta_{N,1}) &\cdots & T_{2N}'(\theta_{N,1})\\
\vdots & & \vdots\\
T_0'(\theta_{N,N}) &\cdots & T_{2N}'(\theta_{N,N})\\
\end{array}\right],
R=\left[\begin{array}{l}0\\ \vdots \\0\\ \lambda \\0\\ \vdots \\
0\end{array}\right].
\]
}
The matrix $T$ is always invertible, which can easily be
deduced from a representation in a Lagrange basis. After
solving (\ref{solvepn}) for a given value of $\lambda$ we
can evaluate
\[
p_N(\pm\tau_i;\ \lambda)=\sum_{i=0}^{2N} \alpha_i
T_i\left(\pm\tau_i\right),\ \ i=0,\ldots,m.
\]
Our implementation is based on a representation in an
orthogonal basis of Chebyshev polynomials.
\begin{remark}\label{remcommon}
We can formally write
\[
p_N(t;\ \lambda)= S(t)(\lambda T-U)^{-1}R,
\]
where $S(t)=[T_0(t)\cdots T_{2N}(t)]$. This shows that
the functions $p_N(\pm\tau_i;\ \lambda)$ are proper
rational function of $\lambda$. Note that the
coefficients of this rational function never need to be
explicitly computed.
\end{remark}

\subsection{Correction of the
$\mathcal{H}_{\infty}$ norm}\label{parcorrect}

We describe how an approximation of
$\|G(j\omega)\|_{\mathcal{H}_{\infty}}$  can be corrected
to the actual value. This is done by solving a set of
nonlinear equations, which we derive first.

\begin{figure}
\begin{center}
\resizebox{10cm}{!}{\includegraphics{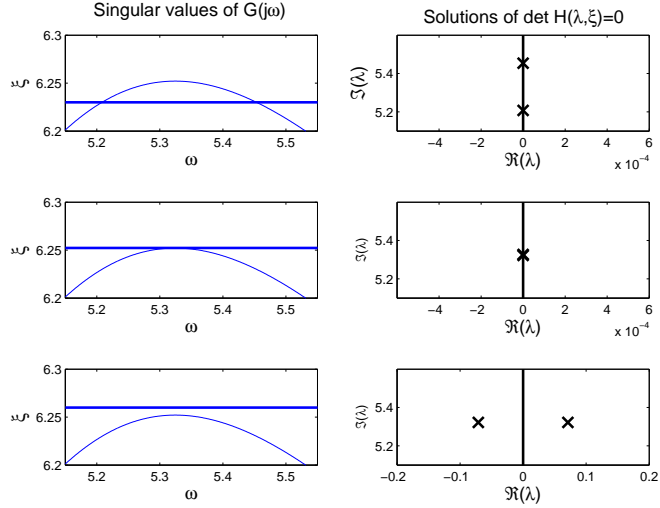}}
 \caption{\label{figD} (left) Intersections of the
 singular value plot of $G$ with the horizontal line $\xi=c$,
 for $c<\hat\xi$ (top),\ $c=\hat\xi$ (middle) and
 $c>\hat\xi$ (bottom). (right) Corresponding eigenvalues of the problem
 (\ref{nonlinear-eigenv2}).}
\end{center}
\end{figure}

Let $\hat\xi\geq 0$ and $\hat\omega\geq 0$ are such that
\begin{equation}\label{direct0}
\|G(j\omega)\|_{\mathcal{H}_{\infty}}=\hat\xi=\sigma_1(G(j\hat\omega)),
\end{equation}
and assume for the moment that the singular value
$\sigma_1(G(j\hat\omega))$ has multiplicity one. Then the
nonlinear eigenvalue problem (\ref{nonlinear-eigenv2}),
with $\xi=\hat\xi$, has a double non-semisimple
eigenvalue $\lambda=j\hat\omega$ (see \cite{lancaster}
for the definition of multiple eigenvalues of a nonlinear
eigenvalue problem). This property is clarified in
Figure~\ref{figD}.
 Therefore, setting
\[
h(\lambda,\xi)=\det H(\lambda,\xi),
\]
the pair $(\hat\omega,\hat\xi)$ satisfies
\begin{equation}\label{direct1}
h(j\omega,\xi)=0,\ \ \frac{\partial
h}{\partial\lambda}(j\omega,\xi)=0.
\end{equation}
%
These complex valued equations seem over-determined but
this is not the case due to the spectral properties of
the operator $\mathcal{L}_{\xi}$. As a corollary of
Proposition~\ref{propsymmetric} we namely get:
\begin{corollary}\label{corextra}
For $\omega\geq 0$, we have
\begin{equation}\label{col1}
\Im\ h(j\omega,\xi)=0
\end{equation}
 and
\begin{equation}\label{col2}
 \Re\ \frac{\partial h}{\partial\lambda}(j\omega,\xi)=0.
\end{equation}
\end{corollary}
\noindent\textbf{Proof.\ } From (\ref{symmetrie}) we get
\[
h(\lambda,\xi)=h(-\lambda,\xi),\ \ \ \frac{\partial
h}{\partial\lambda}( \lambda,\xi)=-\frac{\partial
h}{\partial\lambda} ( -\lambda,\xi).
\]
Substituting $\lambda=j\omega$ yields
\[
\begin{array}{l}
h(j\omega,\xi)=h(-j\omega,\xi)=
\left(h(j\omega,\xi)\right)^*, \\
\frac{\partial h}{\partial\lambda}(j\omega,\xi)=
-\frac{\partial h}{\partial\lambda}(-j\omega,\xi)=
-\left(\frac{\partial
h}{\partial\lambda}(j\omega,\xi)\right)^*
\end{array}
\]
and the assertions follow. \hfill $\Box$

Using Corollary~\ref{corextra} we can simplify the
conditions (\ref{direct1}) to
\begin{equation}\label{direct2x}
\left\{\begin{array}{l}
\Re\ h(j\omega,\xi)=0 \\
\Im\ \frac{\partial h}{\partial\lambda}(j\omega,\xi)=0
\end{array}\right..
\end{equation}
In this way the pair $(\hat\omega,\hat\xi)$ satisfying
(\ref{direct0}) can be directly computed from the two
equations (\ref{direct2x}), for example using Newton's
method, provided that good starting values are available.

The drawback of working directly with (\ref{direct2x}) is
that an explicit expression for the determinant of
$H(\lambda,\xi)$ is required.
To avoid this, let $v_1,v_2\in\mathbb{C}^n$ be such that
\[
H(j\omega,\xi) \left[\begin{array}{c}v_1\\
v_2\end{array}\right]=0,\ \ \ n(v_1,v_2)=0,
\]
where $n(v_1,v_2)=0$ is a normalizing condition. Given
the structure of $H(j\omega,\xi)$ it can be verified that
a corresponding left null vector is given by $[-v_2^*\
v_1^*]$. According to \cite{lancaster} a necessary
condition for the eigenvalue $\lambda=j\omega$ to be
double but non-semisimple is given by\footnote{The
condition (\ref{jordan}) guarantees the existence of a
Jordan chain of length larger than one.}
\begin{equation}\label{jordan}
[-v_2^*\ v_1^*]\ \frac{\partial
H}{\partial\lambda}(j\omega,\xi)
\left[\begin{array}{c}v_1\\
v_2\end{array}\right]=0.
\end{equation}
A simple computation results in
\[
[-v_2^*\ v_1^*]\ \frac{\partial H}{\partial\lambda}
(j\omega,\xi) \left[\begin{array}{c}v_1\\
v_2\end{array}\right]=2\Im\left\{v_2^*\left(I+\sum_{i=1}^p
A_i\tau_i e^{-j\omega\tau_i}+B D_{\xi}^{-1}D^TC\tau_0
e^{j\omega\tau_0}\right)v_1\right\}.
\]
This expression is always real, which is a property
inferred from (\ref{col2}). In this way we end up with
$4n+3$ real equations
\begin{equation}\label{forfinalpref}
\left\{\begin{array}{l}
H(j\omega,\ \xi)\left[\begin{array}{c}v_1\\
v_2\end{array}\right]=0 \\
n(v_1,v_2)=0\\
\Im\left\{v_2^*\left(I+\sum_{i=1}^p A_i\tau_i
e^{-j\omega\tau_i}+B D_{\xi}^{-1}D^TC\tau_0
e^{j\omega\tau_0}\right)v_1\right\}=0
\end{array}\right.
\end{equation}
in  the $4n+2$ unknowns $\Re(v_1),\Im
(v_1),\Re(v_2),\Im(v_2),\omega$ and $\xi$. These
equations are still overdetermined because  the property
(\ref{col1}) is not explicitly exploited in the
formulation, unlike the property (\ref{col2}). However,
this property makes the equations (\ref{forfinalpref})
exactly solvable.

In our implementation we solve the equation
(\ref{forfinalpref}) in least squares sense using the
Gauss Newton algorithm, which can be shown to be
quadratically converging in the case under consideration
where the residual in the desired solution is zero. The
program gives a warning when the correction involves a
relative change larger than 10\%, because this
 indicates that the approximation in the prediction step
 might not be accurate enough.


%
%

\medskip

The above results can be used to compute the
$\mathcal{H}_{\infty}$ norm of $G(j\omega)$, in the
following way. Given an approximation of a pair
$(\hat\xi,\hat\omega)$ satisfying (\ref{direct0}) and
given corresponding estimates of $v_1$ and $v_2$, we can
find the exact values by solving the nonlinear equations
(\ref{forfinalpref}) in least squares sense. The fact
that the residual must be zero in the desired solution
can be used an additional optimality certificate. The
approximation $(\hat\xi,\hat\omega)$ can be obtained by
the prediction step outlined in $\S\ref{parpredict}$.
When using Algorithm~\ref{algmiddle} for the prediction
step the total algorithm becomes:

\begin{algorithm}\label{algtotal} \ \

{\small
 {\em
 \noindent Input: system data, $N$, symmetric
 grid $\Omega_N$ defined by (\ref{chebpoints}), candidate
critical frequency $\omega_t$\\
\hspace*{1cm}if available, tolerance tol for prediction step\\
 Output: $\|G(j\omega)\|_{\mathcal{H}_{\infty}}$
 }
\begin{enumerate}
\item[]\hspace*{-0.6cm} \underline{\emph{Prediction step:}}
\item[]\hspace*{-0.6cm} Apply Algorithm~\ref{algmiddle}.
\item[]\hspace*{-0.6cm} \underline{\emph{Correction step:}}
\item determine all eigenvalues
$\{j\omega^{(1)},\ldots,j\omega^{(p)}\}$ of
 $\mathcal{L}_{\xi_l}^N$ on the positive imaginary axis,
 \item for all $i\in\{1,\ldots,p\}$, solve (\ref{forfinalpref}) with
 starting values
 \[
\omega=\omega^{(i)},\ \ \xi=\xi_l,\
\left[\begin{array}{c}v_1\\
v_2\end{array}\right]=
\arg\min{\|H(j\omega^{(i)},\xi_l)v\|}/{\|v\|};
 \]
denote the solution with $(\tilde u^{(i)},\tilde
 v^{(i)},\tilde\omega^{(i)},\tilde\xi^{(i)})$
\item set
$\|G(j\omega)\|_{\mathcal{H}_{\infty}}:=\max_{1\leq i\leq
p}\tilde\xi^{(i)}$
\end{enumerate}
}
\end{algorithm}

\smallskip


\subsection{Number of discretization
points}\label{parsteplength}

An important aspect in the application of
Algorithm~\ref{algtotal} consists of choosing $N$, or,
equivalently, the number of grid points for the
discretization of $\mathcal{L}_{\xi}$ into
$\mathcal{L}_{\xi}^N$. On the one hand, from a
computational point of view, $N$ should be as small as
possible, given that the computational cost of
Algorithm~\ref{algtotal} is dominated by the computation
of the eigenvalues of the matrix $\mathcal{L}_{\xi}^N$ with dimensions ($2N+1)2n$.
On the other hand, $N$ should be sufficiently large to
generate starting values for which the corrector
converges to the desired values.
In other but similar problems, in particular the problem
of the computation of characteristic roots, mostly a
heuristic or a 'safe' overestimate for the number of
discretization points is used.
%
%
In our case, the relations established in
Sections~\ref{sectheo}-\ref{par1} and summarized in
Table~\ref{tabelleke} turn out to be very useful in the
determination of $N$. This is explained in what follows.

\begin{figure}
\begin{center}
\resizebox{12cm}{!}{\includegraphics{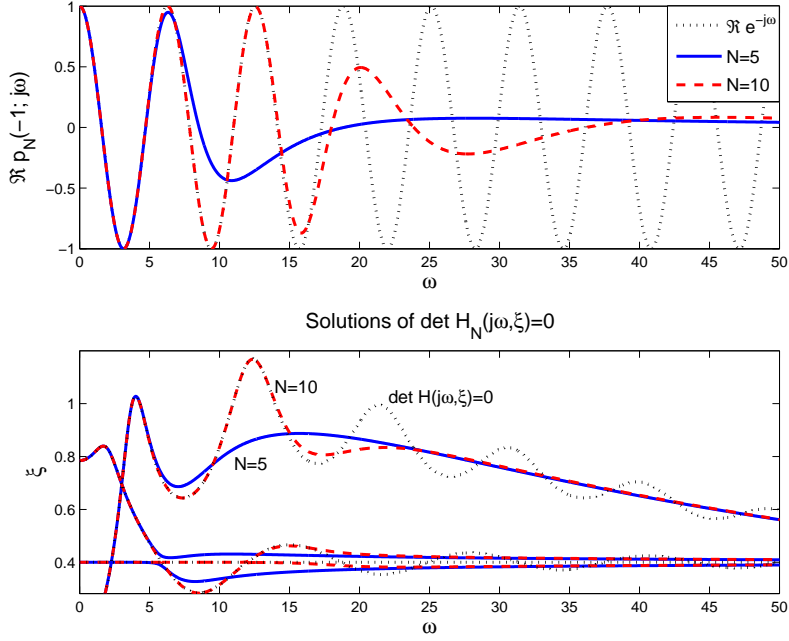}}
 \caption{\label{figE} (top) Comparison between the
 functions
 $\omega\mapsto\exp(j\omega)$ and $\omega\mapsto
 p_N(-1;\ j\omega)$ for $N=5$ and $N=10$. (bottom)
 Corresponding comparison between the solutions of $\det
 H(j\omega,\xi)=0$ (the singular value plot of $G$) and
 the solutions of $\det H_N(j\omega;\ \xi)=0$, with problem
 data (\ref{example1}). An approximation with $N=10$ is
 sufficient to cover the highest peak, unlike an
 approximation with $N=5$.
 }
\end{center}
\end{figure}

From Definition~\ref{defpn} the function
\begin{equation}\label{fucpn}
t\in[-1,\ 1]\mapsto p_N(t;\  \lambda)
\end{equation}
is an approximation of the function
\begin{equation}\label{fucpn2}
t\in[-1,\ 1]\mapsto e^{\lambda t}.
\end{equation}
Moreover, a comparison between $H(\xi,\lambda)$, defined
by (\ref{defh}), and $H_N(\xi,\lambda)$, defined by
(\ref{defhn}), learns that the effect of approximating
$\mathcal{L}_{\xi}$ with $\mathcal{L}_{\xi}^N$ can be
interpreted as the effect of approximating the
exponential functions $e^{\pm\lambda\tau_i}$ with the
rational functions $p_N(\pm\tau_i;\ \lambda)$ for
$i=0,\ldots,m$. Hence, the number $N$ in the prediction
step of Algorithm~\ref{algtotal} should be chosen in such
a way that
\[
p_N(\pm\tau_i;\ j\omega)\approx e^{\pm j\omega\tau_i},\
i=0,\ldots,m,
\]
in the relevant frequency range, that is, where the
highest peak values in the singular value plot of
$G(j\omega)$ occur. This is illustrated in
Figure~\ref{figE}.

\begin{figure}
\begin{center}
\resizebox{7.5cm}{!}{\includegraphics{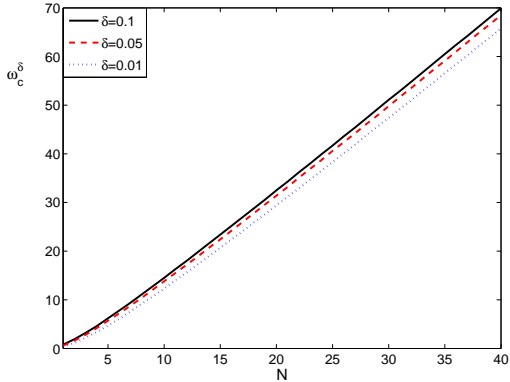}}
 \caption{\label{figF} The cut-off frequency $\omega_c^{\delta}$, as a function of $N$, for
 $\delta=0.01$, $\delta=0.05$ and $\delta=0.10$.}
\end{center}
\end{figure}

In Figure~\ref{figF} we depict the 'cut-off
frequency'
\begin{equation}\label{cut-off}
\omega_c^{\delta}(N):=\min\left\{\omega\geq 0:\
\max_{t\in[-1,\ 1]}|e^{j\omega t}-p_N(t;\
j\omega)|\geq\delta\right\}
\end{equation}
as a function of $N$ for different values of $\delta$.
The importance is as follows: if $\omega\in[0,\
\omega_c^{\delta}(N)]$, then it is guaranteed that
$p_N(\pm\tau_i; j\omega)$ approximates $e^{\pm
j\omega\tau_i}$ with relative error smaller than
$\delta$, for $i=0,\ldots,m$.
 Because $\omega_c^{\delta}$ is \emph{independent} of the problem,
 Figure~\ref{figF}
can assist in determining a suitable value of $N$, when
information about the relevant frequency range in the
singular value plot of $G$ is available. It is not
necessary to choose $\delta$ very small, because the
prediction step in Algorithm~\ref{algmiddle} is followed
by a correction step, hence, it is in the first place
sufficient to capture qualitatively the highest peaks in
the singular value plot.
%
%
In our implementation of Algorithm~\ref{algtotal} the
value $\delta=0.1$ is taken. The user has the option to
specify the corresponding cut-off frequency
$\omega_c^{0.1}$, from which the value of $N$ is
automatically determined. Otherwise the default value
$N=15$ is used. If in the prediction step a peak value is
computed at a frequency $\omega_p$ satisfying
$\omega_p>\omega_c^{0.1}$, then a warning is generated
with an advice to increase the cut-off frequency.


\smallskip

Our experience from extensive benchmarking learns that in
most practical problems a very small value of $N$ can be
taken (the default $N=15$ is largely sufficient). This
can be motivated as follows. First, via Corollary
\ref{colcompute} the $\mathcal{H}_{\infty}$ problem is
transformed into a problem of checking the eigenvalues of
$\mathcal{L}_{\xi}$ on the imaginary axis. Because these
eigenvalues  are typically among the \emph{smallest}
eigenvalues of $\mathcal{L}_{\xi}$ and the individual
eigenvalues of $\mathcal{L}_\xi^N$ exhibit \emph{spectral
convergence}
 to the corresponding eigenvalues of $\mathcal{L}_{\xi}$
(this can be shown following the lines of
\cite{breda:nonlocal}), a small value of $N$ already
yields good approximations of the imaginary axis
eigenvalues.
Second, in most applications the delay parameters are
critical from a stability point of view, in the sense
that increasing the delay parameters eventually
destabilizes the system, but stability is maintained by
decreasing the delay parameters (often referred to in the
literature with the term delay-dependent stability). The
connection between delay-dependent stability and the
possibility to work with a small value of $N$ in the
$\mathcal{H}_{\infty}$ computation is intuitively
explained with an example.
\begin{example}\label{excritical}
We consider the transfer function
\[
G(j\omega)=C(j\omega I-A_0-A_1 e^{-j\omega\tau})^{-1}B+D,
\]
with nominal delay $\tau=\tau_n=1$. We assume that the
system
\begin{equation}\label{defGex}
\dot{x}(t)=A_0 x(t)+A_1 x(t-\tau)
\end{equation}
is stable for all $\tau\in [0,\ 1]$ but not
delay-independent stable. It is well known that, if the
delay is increased, characteristic roots can only cross
the imaginary axis at a finite number of points, say
$\{j\Omega_1,j\Omega_2,\ldots\}$ (see e.g.\
\cite{bookwim}). The periodicity of the presence
of an eigenvalue at $j\Omega_i$ with respect to delay
shifts of $2\pi/\Omega_i$, implies
\begin{equation}\label{estau0}
\tau_n=1\leq\frac{2\pi}{\max_i\Omega_i},
\end{equation}
which leads to
\[
1\geq\frac{\max_i\Omega_i}{2\pi}.
\]
If we choose $N=8$, then we find from Figure~\ref{figF}
that $\omega_c^{0.1}>10$. As a consequence, we get
\[
\omega_c^{0.1}>10\geq
\frac{10}{2\pi}\max_i\Omega_i\approx 1.59\
\max_i\Omega_i.
\]
This is expected to be sufficient for a correct
computation of the $\mathcal{H}_{\infty}$ norm of $G$
because the highest peak values in the singular value
plot of (\ref{defGex}) typically occur in the frequency
range where characteristic roots can come close to or
cross the imaginary axis. Moreover, the factor $1.59$ does
not take into account the conservatism of the estimate
(\ref{estau0}), which can be very large.
\end{example}

\section{Numerical examples}\label{secex}

Throughout the paper, the obtained results and the ideas
behind Algorithm~\ref{algtotal} have been illustrated
with the problem data
\begin{eqnarray} \label{example1}
\nonumber A_0&=&\left[
  \begin{array}{rrr}
    108 & 110 & 18 \\
    -107 & -109 & -17 \\
    -217 & -217 & -37 \\
  \end{array}
\right],
 A_1=\left[
  \begin{array}{rrr}
    46.5 & 46.5 & 1.5 \\
    -46.5 & -46.5 & -1.5 \\
    -93 & -93 & -3 \\
  \end{array}
\right], \\
\nonumber A_2&=&\left[
  \begin{array}{rrr}
    -0.3 & 0.3 & -0.3 \\
    0.3 & -0.3 & 0.3 \\
    0 & 0 & 0 \\
  \end{array}
\right], B=\left[
    \begin{array}{rr}
      0.5 & -90 \\
      -0.5 & 90 \\
      0 & 180 \\
    \end{array}
  \right], \\
\nonumber C&=&\left[
    \begin{array}{rrr}
      1 & -1 & 1 \\
      0.1833 & 0.1833 & 0.1833 \\
    \end{array}
  \right],
D=\left[
               \begin{array}{rr}
                 0.4 & 0 \\
                 0 & 0 \\
               \end{array}
             \right], \\
\nonumber \tau_0&=&0.5,\quad \tau_1=0.667,\quad \tau_2=1.
\end{eqnarray}
In particular, the singular value plot of the
corresponding transfer function $G$ is shown in
Figure~\ref{figA} and the solutions of $\det
H_N(j\omega,\xi)=0$ in Figure~\ref{figB}, for $N=10$. The
two parameter search in the prediction step is
illustrated in Figure~\ref{figC}. The effect of the
approximation in the prediction step is shown in
Figure~\ref{figE}.  With $N=10$, the predicted $\Hi$ norm
is $\xi_{pred}=1.1626$ and the corrected $\Hi$ norm is
given by $\xi_{corr}=1.1696$.

In Table~\ref{table:Benchmarks} we present the results of
benchmarking of our code with $12$ problems.  The second
column shows the size of matrices $A_i$, $n$, and the
number of state delays, $m$. The third column gives the
minimum value of $N$ such that in the correction step the
desired solution is computed. The fourth and fifth
columns contain the predicted and corrected $\Hi$ norms
of the corresponding time-delay system. The plant $G12$
correspond to the problem data (\ref{example1}).

\begin{table}[h!]
\begin{center}
\begin{tabular}{l|l|l|l|l}
  Plants & $(n,m)$ & $N$ & $\xi_{pred}$ & $\xi_{corr}$ ($\|G(j\omega)\|_{\mathcal{H}_{\infty}}$) \\
  \hline
  \hline
  $G1$ & $(3,1)$ & $1$ & $10.0235$ & $10.0235$ \\
  $G2$ & $(3,1)$ & $3$ & $3.3693$ & $3.3709$ \\
  $G3$ & $(1,1)$ & $2$ & $0.7158$ & $0.7196$ \\
  $G4$ & $(1,1)$ & $3$ & $1.9230$ & $1.9883$ \\
  $G5$ & $(3,3)$ & $2$ & $0.8852$ & $0.8848$ \\
  $G6$ & $(3,3)$ & $5$ & $0.8974$ & $0.9356$ \\
  $G7$ & $(4,3)$ & $10$ & $1.6259$ & $1.6283$ \\
  $G8$ & $(10,7)$ & $4$ & $22.2979$ & $22.3195$ \\
  $G9$ & $(20,9)$ & $13$ & $1.2827$ & $1.2903$ \\
  $G10$ & $(40,3)$ & $3$ & $811.0898$ & $814.6221$ \\
  $G11$ & $(3,2)$ & $9$ & $1.1579$ & $1.1696$ \\
  $G12^{*}$ & $(3,2)$ & $18$ & $1.1626$ & $1.1696$ \\
  \hline
\end{tabular}
\end{center}
\caption{Benchmarks for the $\Hi$ norm computation.}
\label{table:Benchmarks}
\end{table}

For the plant $G12$ a warning is generating when using
the default value $N=15$, indicating that the
corresponding cut-off frequency might be too small. With
$N=18$ no warning is given and the results are correct.
We note that this example has been constructed in such a
way that a relatively large value of $N$ is necessary, by
taking very large values of 'non-critical' delay
parameters (see the discussion on critical delay
parameters in the paragraph before
Example~\ref{excritical}). When the delay parameters are
critical from a stability point of view, which is the
case in most practical problems, a much smaller values of
$N$ is sufficient, as motivated with
Example~\ref{excritical}.

\smallskip

The problem data for the above benchmark examples and a
MATLAB implementation of our code for the $\Hi$ norm
computation are available at the website
\begin{verbatim}
http://www.cs.kuleuven.be/~wimm/software/hinf/
\end{verbatim}

\section{Discussion of alternative approaches and concluding remarks}\label{seccon}

In this article we have described an algorithm for the
computation of $\Hi$ norms for time-delay systems, which
relies on a two-step approach: a prediction step where an
approximation is computed based on a finite-dimensional
approximation, and a local corrector. It should be
noticed that these two steps are to some extent
independent of each other. In particular, other choices
for a finite-dimensional approximation in the prediction
step are possible. While our approach is based on a
discretization of $\mathcal{L}_{\xi}$ it is for instance
also possible to use a direct approximation of the
transfer function $G$, by replacing the exponential
functions with a rational approximation. In this context
our choice for a spectral discretization of
$\mathcal{L}_{\xi}$ is motivated as follows.

\begin{itemize}
\item An approach based on a spectral discretization of an appropriately
defined derivative operator with nonlocal boundary
conditions is known to be not only an accurate
(cf.~spectral convergence of the eigenvalues) but also a
numerically stable way to solve infinite-dimensional
eigenvalue problems in the context of time-delay system
(see \cite{breda:nonlocal}, see also \cite{phdverheyden}
for a discussion on various methods for the "dual"
problem of computing characteristic roots). On the
contrary, working with an explicit rational approximation
of $G$ may lead to an ill-conditioned Hamiltonian matrix
when applying Proposition~\ref{propintro} to the
resulting finite-dimensional system. This is due to
potential large differences in magnitudes of the
coefficients in rational approximants of high order (a
high order is necessary for "globally" capturing the
transfer function in the relevant frequency range).
\item From Definition~\ref{Defpn} and Property~\ref{propnl}
the effect of discretizing $\mathcal{L}_{\xi}$ can be
interpreted as the effect of approximating
$e^{\pm\lambda\tau_i}$ in $H(\lambda,\xi)$ by
$p_N(\pm\tau_i;\ \lambda)$, where the function
(\ref{fucpn})
is obtained as a polynomial approximation of
(\ref{fucpn2}),
satisfying collocation and interpolation conditions on
the grid $\Omega_N$,\ i.e.\ the exponential function is
approximated over the full interval $[-1\ 1]$, to which
all delays belong. As a consequence, the dimensions of
the matrix $\mathcal{L}_{\xi}^N$, $(2N+1)2n\times
(2N+1)2n$, are \emph{independent} of the number of delays
in the problem, $m$. This can also be seen from
Remark~\ref{remcommon}, which shows that the poles of the
rational functions $\lambda\mapsto p_N(-\tau_i;\
\lambda)$ are independent of $\tau_i$, that is, one can
interpret the effect of discretizing $\mathcal{L}_{\xi}$
as the effect of an approximation of all exponential
functions by rational functions with common poles. With
Pad\'e and many other types of rational approximations
these poles are not the same, and the dimension of the
discretized systems will become proportional to the
number of delays.
\item An important advantage of a
direct rational approximation of the exponential
functions in (\ref{transfer}) is that a high accuracy in
a relevant frequency range can easily be guaranteed by
the choice of the order of the approximation. We have
demonstrated that this is also possible when working with
a spectral discretization of the operator
$\mathcal{L}_{\xi}$, once again via  the interpretation
of the effect of its discretization as the effect of a
rational approximation (although the coefficients of the
rational functions $p_N(-\tau_i;\ \lambda)$ never needed
to be explicitly computed). This property was used in
\S\ref{parsteplength} for the determination of the number
of discretization points.
\end{itemize}

\smallskip

The algorithm has been intensively tested, and turns out
to be very robust. The computational cost is dominated by
the determination of the eigenvalues of the matrix
$\mathcal{L}_{\xi}$. In our current implementation all
eigenvalues are computed. However, since the algorithm
only needs the eigenvalues of $\mathcal{L}_{\xi}$ in the
vicinity of the imaginary axis, which are typically among
the smallest eigenvalues, subspace methods based on
inverse iteration become appealing for large problems.
Instrumental to this the techniques described in
Section~2.2 of \cite{phdverheyden} allow to bring the
eigenvalue problem of $\mathcal{L}_{\xi}$  in a form for
which matrix vector become cheap. These issues, as well
as the application of the algorithm to
$\mathcal{H}_{\infty}$ synthesis problems, are outside
the scope of this paper.

\section*{Acknowledgements}
This article present results of the Belgian Programme on
Interuniversity Poles of Attraction, initiated by the
Belgian State, Prime Minister's Office for Science,
Technology and Culture, and of the Optimization in
Engineering Centre OPTEC. The authors wish to thank the editors and anonymous reviewers for their careful reading and their constructive comments to improve the quality and readability of the paper.

\bibliographystyle{plain}
\bibliography{siampaper}

\end{document}